\documentclass[10pt,onecolumn]{article}

\usepackage{times}
\usepackage{latexsym}

\usepackage{amsmath} 
\usepackage{amssymb} 

\makeatletter
\def\@maketitle{%
  \vbox to 6.5cm{%
    \hsize\textwidth
    \linewidth\hsize
    \vspace{1.5cm}
    \centering
    {\bfseries\LARGE \@title \par}
    \vspace{12pt}
    {\fontsize{11pt}{13pt}\selectfont \begin{tabular}[t]{c}\@author \end{tabular}\par}
    \vfill} 
}
\renewcommand\section{\@startsection{section}{1}{\z@}%
                       {-12\p@ \@plus -4\p@ \@minus -4\p@}%
                       {6\p@ \@plus 4\p@ \@minus 4\p@}%
                       {\normalfont\large\bfseries
                        \rightskip=\z@ \@plus 8em\pretolerance=10000 }}
\renewcommand\subsection{\@startsection{subsection}{2}{\z@}%
                       {-12\p@ \@plus -4\p@ \@minus -4\p@}%
                       {6\p@ \@plus 4\p@ \@minus 4\p@}%
                       {\normalfont\fontsize{11pt}{13pt}\selectfont\bfseries
                        \rightskip=\z@ \@plus 8em\pretolerance=10000 }}
\renewcommand\subsubsection{\@startsection{subsubsection}{3}{\z@}%
                       {-12\p@ \@plus -4\p@ \@minus -4\p@}%
                       {6\p@ \@plus 4\p@ \@minus 4\p@}%
                       {\normalfont\normalsize\itshape}}
\renewcommand\paragraph{\@startsection{paragraph}{4}{\z@}%
                       {-12\p@ \@plus -4\p@ \@minus -4\p@}%
                       {-0.5em \@plus -0.22em \@minus -0.1em}%
                       {\normalfont\normalsize\itshape}}
\makeatother

\renewenvironment{abstract}%
  {\small
    \list{}{\labelwidth0pt
      \leftmargin0pt \rightmargin\leftmargin
      \listparindent\parindent \itemindent0pt
      \parsep0pt
      }%
    \item[\hskip\labelsep\bfseries\abstractname\enspace --] \itshape}{\endlist}

\newcommand{\keywordsname}{Keywords}
\newenvironment{keywords}%
  {\small
    \list{}{\labelwidth0pt
      \leftmargin0pt \rightmargin\leftmargin
      \listparindent\parindent \itemindent0pt
      \parsep0pt
      }%
    \item[\hskip\labelsep\bfseries\keywordsname:]}{\endlist}

\begin{document}

\pagestyle{plain} 

\title{On the Tweety Penguin Triangle Problem}

\author{\begin{tabular}{c@{\extracolsep{8em}}c}
{\bf Jean Dezert} & {\bf Florentin Smarandache}\\
ONERA & Department of Mathematics\\
29 Av. de la  Division Leclerc & University of New Mexico\\
92320 Ch\^{a}tillon & Gallup, NM 8730\\
France & U.S.A.\\
{\tt Jean.Dezert@onera.fr} & {\tt smarand@unm.edu}
\end{tabular}}

\date{}
\maketitle
\pagestyle{plain}


\begin{abstract}
In this paper, one studies the famous well-known and challenging Tweety Penguin Triangle Problem (TPTP or TP2) pointed out by Judea Pearl in one of his books. We first present the solution of the TP2  based on the fallacious  Bayesian reasoning and prove that  reasoning cannot be used to conclude on the ability of the penguin-bird Tweety to fly or not to fly. Then we present in details the counter-intuitive solution obtained from the Dempster-Shafer Theory (DST). Finally, we show how the solution can be obtained  with our new theory of plausible and paradoxical reasoning (DSmT).
\end{abstract}



\section{Introduction}

Judea Pearl claimed that DST of evidence fails to provide a reasonable solution for the combination of evidence even for apparently very simple fusion problem \cite{Pearl_1988,Pearl_1990}. Most criticisms are answered by Philippe Smets in \cite{Smets_1991,Smets_1992}. 
The Tweety Penguin Triangle Problem (TP2) is one of the typical exciting and challenging problem for all theories managing uncertainty and conflict because it shows the real difficulty to maintain truth for automatic reasoning systems when the classical property of transitivity (which is basic to the material-implication) does not hold. In his book \cite{Pearl_1988}, Judea Pearl presents and discusses in details the semantic clash between Bayes vs. Dempster-Shafer reasoning. We present here our new analysis on this problem and provide a solution of the Tweety Penguin Triangle Problem based on our new theory of plausible and paradoxical reasoning, known as DSmT (Dezert-Smarandache Theory). We show how this problem can be attacked and solved by our new reasoning with help of the (hybrid) DSm rule of combination \cite{Dezert_2004Book}.\\

The purpose of this paper is not to browse all approaches available in literature for attacking the TP2 problem but only to provide a comparison of the DSm reasoning with respect to the  Bayesian reasoning and to the plausible reasoning of DST framework. Interesting but complex analysis on this problem based on default reasoning and $\epsilon$-belief functions can be also found by example in \cite{Smets_1991} and \cite{Benferhat_2000}. Other interesting and promising issues for the TP2 problem based on the fuzzy logic of Zadeh \cite{Zadeh_1978} jointly with the theory of possibilities \cite{Dubois_1980,Dubois_1987} are under investigations. Some theoretical research works on new conditional event algebras (CEA) have emerged in literature \cite{Goodman_1997} since last years and could offer a new track  for attacking the TP2 problem although unfortunately no clear didactic, simple and convincing examples are provided to show the real efficiency and usefulness of these theoretical investigations.

\section{The Tweety Penguin Triangle Problem }

This very important and challenging problem, as known as the Tweety Penguin Triangle Problem (TP2) in literature, is presented in details by Judea Pearl in \cite{Pearl_1988}. We briefly present here the TP2 and the solutions based first on fallacious  Bayesian reasoning and then on the Dempster-Shafer reasoning. We will then focus our analysis of this problem from the DSmT framework and the DSm reasoning. \\

\noindent
Let's consider the set $R=\{r_1,r_2,r_3\}$ of given rules:
\begin{itemize}
\item $r_1$: "{\bf{Penguins normally don't fly}}" $\Leftrightarrow$ $(p\rightarrow \neg f)$
\item $r_2$: "{\bf{Birds normally fly}}" $\Leftrightarrow$ $(b\rightarrow f)$
\item $r_3$: "{\bf{Penguins are birds}}" $\Leftrightarrow$ $(p \rightarrow b)$
\end{itemize}

\noindent
To emphasize our strong conviction in these rules we commit them some high confidence weights $w_1$, $w_2$ and $w_3$ in $[0,1]$ with 
$w_1=1-\epsilon_1$, $w_2=1-\epsilon_2$ and $w_3=1$ (where $\epsilon_1$ and $\epsilon_2$ are small positive quantities). The conviction in these rules is then represented by the set $W=\{w_1,w_2,w_3\}$ in the sequel.\\

Another useful and general notation adopted by Judea Pearl in the first pages of his book \cite{Pearl_1988} to characterize these three weighted rules is the following one (where $w_1,w_2,w_3 \in [0,1]$):
$$r_1: p\overset{w_1}{\rightarrow}(\neg f) \qquad r_2: b\overset{w_2}\rightarrow f 
 \qquad r_3: p\overset{w_3}\rightarrow b$$

When $w_1,w_2,w_3 \in \{0,1\}$ the classical logic is the perfect tool to conclude on the truth or on the falsity of a proposition built from these rules  based on the standard propositional calculus mainly with its three fundamental rules (Modus Ponens, Modus Tollens and Modus Barbara - i.e. transitivity rule).
When $0<w_1,w_2,w_3<1$, the classical logic can't be applied because the Modus Ponens, the Modus Tollens and the Modus Barbara do not longer hold and some other tools must be chosen. This will discussed in detail in section \ref{SEC:PEARLWEAKNESS}.\\

\noindent{\bf{Question:}} Assume we observe an animal called Tweety (T) that is categorically classified as a bird (b) and a penguin (p), i.e. our observation is $O\triangleq [T=(b\cap p)]=[(T=b)\cap (T=p)]$. The notation $T=(b\cap p)$ stands here for "Entity $T$ holds property $(b\cap p)$". What is the belief (or the probability - if such probability exists) that Tweety can fly given the observation $O$ and all information available in our knowledge base (i.e. our rule-based system $R$ and $W$) ?\\

The difficulty of this problem for most of artificial reasoning systems (ARS) comes from the fact that, in this example, the property of transitivity, usually supposed satisfied from material-implication interpretation \cite{Pearl_1988}, $(p \rightarrow b,b\rightarrow f) \Rightarrow (p \rightarrow f)$ does not hold here (see section \ref{SEC:PEARLWEAKNESS}). In this interesting example, the classical property of inheritance is thus broken. Nevertheless a powerful artificial reasoning system must be able to deal with such kind of difficult problem and must provide a reliable conclusion by a general mechanism of reasoning whatever the values of convictions are (not only restricted to values close to either 0 or 1). We examine now three ARS based on the Bayesian reasoning \cite{Pearl_1988} which turns to be fallacious and actually not appropriate for this problem and we explain why, on the Dempster-Shafer Theory (DST) \cite{Shafer_1976} and on the Dezert-Smarandache Theory (DSmT) \cite{Dezert_2004Book}.

\section{The fallacious Bayesian reasoning}

We first present the fallacious Bayesian reasoning solution drawn from the J. Pearl's book in \cite{Pearl_1988} (pages 447-449) and then we explain why the solution which seems at the first glance correct with intuition is really fallacious. We then explain why the common rational intuition turns actually to be wrong.

\subsection{The Pearl's analysis}

To preserve mathematical rigor, we introduce explicitly all information available in the derivations.
 In other words, one wants to evaluate using the Bayesian reasoning, the conditional probability, if it exists, $P(T=f | O, R, W)=P(T=f | T=p,T=b,R,W)$. The Pearl's analysis is based on the assumption that a conviction on a given rule can be interpreted as a conditional probability (see \cite{Pearl_1988} page 4). In other words if one has a given rule $a \overset{w}{\rightarrow} b$ with $w\in [0,1]$ then one can interpret, at least for the calculus, $w$ as $P(b|a)$ and thus the probability theory and Bayesian reasoning can help to answer to the question. We prove in the following section that such model cannot be reasonably adopted. For now, we just assume that such probabilistic model holds effectively as Judea Pearl does. Based on this assumption, since the conditional term/information $(T=p,T=b,R,W)$ is strictly equivalent to $(T=p,R,W)$ because of the knowledge of rule $r_3$ with certainty (since $w_3=1$), one gets easily the fallacious intuitive expected Pearl's result:
 \begin{align*}
P(T=f | O,R,W) &= P(T=f | T=p,T=b,R,W)\\
P(T=f | O,R,W) &\equiv P(T=f | T=p,R,W) \\
P(T=f | O,R,W) & =1-P(T=\neg f | T=p,R,W) \\
P(T=f | O,R,W) & = 1 - w_1=\epsilon_1
 \end{align*}
From this simple analysis, the Tweety's "birdness" does not render her a better flyer than an ordinary penguin as intuitively expected and the probability that Tweety can fly remains very low which looks  normal. We reemphasize here the fact, that in his Bayesian reasoning J. Pearl assumes that the weight $w_1$ for the conviction in rule $r_1$ can be interpreted in term of a real probability measure $P( \neg f | p)$.  This assumption is necessary to provide the rigorous derivation of $P(T=f | O,R,W)$. It turns out however that convictions $w_i$ on logical rules cannot be interpreted in terms of probabilities as we will prove in the next section.\\

When rule $r_3$ is not asserted with absolute certainty (i.e. $w_3=1$) but is subject to exceptions, i.e. $w_3=1-\epsilon_3 < 1$, the fallacious Bayesian reasoning yields (where notations $T=f$, $T=b$ and $T=p$ are replaced by $f$, $b$ and $p$ for notation convenience):
\begin{align*}
P(f | O,R,W) &= P(f | p,b,R,W)\\
P(f | O,R,W)&= \frac{P(f,p,b | R,W)}{P(p,b | R,W)}\\
P(f | O,R,W)&= \frac{P(f,b | p,R,W)P(p |R,W)}{P(b | p,R,W)P(p|R,W)}
\end{align*}
\noindent By assuming $P(p |R,W)> 0$, one gets after simplification by $P(p |R,W)$
\begin{align*}
P(f | O,R,W) &= \frac{P(f,b | p,R,W)}{P(b | p,R,W)}\\
P(f | O,R,W)&= \frac{P(b | f,p,R,W)P(f | p,R,W)}{P(b | p,R,W)}
\end{align*}
If one assumes $P(b | p,R,W)=w_3=1-\epsilon_3$ and $P(f | p,R,W)= 1-P(\neg f | p,R,W)=1- w_1=\epsilon_1$, one gets
$$P(f | O,R,W) = P(b | f,p,R,W)\times \frac{\epsilon_1}{1-\epsilon_3}$$
\noindent
Because $0\leq P(b | f,p,R,W)\leq 1$, one finally gets the Pearl's result \cite{Pearl_1988} (p.448)
\begin{equation}
P(f | O,R,W) \leq  \frac{\epsilon_1}{1-\epsilon_3}
\end{equation}

\noindent which states that the observed animal Tweety (a penguin-bird) has a very small probability of flying as long as $\epsilon_3$ remains small, regardless of how many birds cannot fly ($\epsilon_2$), and has consequently a high probability of not flying because $P(f | O,R,W)+P(\bar{f} | O,R,W)=1$ since the events $f$ and $\bar{f}$ are mutually exclusive and exhaustive (assuming that the Pearl's probabilistic model holds ... ).

\subsection{The weakness of the Pearl's analysis}
\label{SEC:PEARLWEAKNESS}
We prove now that the previous Bayesian reasoning is really fallacious and the problem is truly undecidable to conclude about the ability of Tweety to fly or not to fly if a deep analysis is done. Actually,  the Bayes' inference is not a classical inference \cite{Dambreville_2004}. Indeed, before applying blindly the Bayesian reasoning as in the previous section, one first has to check that the probabilistic model is well-founded to characterize the convictions of the rules of the rule-based system under analysis. We prove here that such probabilistic model doesn't hold for a suitable and useful representation of the problem and consequently for any problems based on the weighting of logical rules (with positive weighting factors/convictions below than 1).

\subsubsection{Preliminaries}

We just remind here only few important principles of the propositional calculus of the classical Mathematical Logic which will be used in our demonstration. 
A simple notation, which may appear as unusual for logicians, is adopted here just for convenience. 
 A detailed presentation of the propositional calculus and Mathematical Logic can be easily found in many standard mathematical textbooks like \cite{Nidditch_1962,Mendelson_1997,Kleene_2002}. Here are these important principles:
\begin{itemize}
\item {\bf{Third middle excluded principle}} :  A logical variable is either true or false, i.e.
\begin{equation}
a \vee \neg a
\label{eq:TEM}
\end{equation}
\item {\bf{Non-contradiction law}} :  A logical variable can't be both true and false, i.e.
\begin{equation}
\neg (a \wedge \neg a)
\label{eq:NCL}
\end{equation}
\item {\bf{Modus Ponens}} :  This rule of the propositional calculus states that if a logical variable $a$ is true and  $a\rightarrow b$ is true, then $b$ is true (syllogism principle), i.e.
\begin{equation}
(a \wedge (a \rightarrow b) ) \rightarrow b
\label{eq:MP}
\end{equation}
\item {\bf{Modus Tollens}} :  This rule of the propositional calculus states that if a logical variable $\neg b$ is true and  $a\rightarrow b$ is true, then $\neg a$ is true, i.e.
\begin{equation}
(\neg b \wedge (a \rightarrow b) ) \rightarrow \neg a
\label{eq:MT}
\end{equation}
\item {\bf{Modus Barbara}} :  This rule of the propositional calculus states that if $a\rightarrow b$ is true and $b\rightarrow c$ is true then $a\rightarrow c$ is true (transitivity property), i.e.
\begin{equation}
((a \rightarrow b) \wedge (b \rightarrow c) ) \rightarrow (a \rightarrow c)
\label{eq:MB}
\end{equation}
\end{itemize}

From these principles, one can prove easily, based on the truth table method, the following property (more general deducibility theorems in Mathematical Logic can be found in \cite{Smarandache_1979,Smarandache_1995}) :
\begin{equation}
((a \rightarrow b) \wedge (c \rightarrow d)) \rightarrow ((a\wedge c)\rightarrow (b\wedge d))
\label{eq:MI}
\end{equation}

\subsubsection{Analysis of the problem when   $\epsilon_1=\epsilon_2=\epsilon_3=0$}

We first examine the TP2 when one has no doubt in the rules of our given rule-based systems, i.e.
\begin{equation*}
\begin{cases}
r_1: p\overset{w_1=1-\epsilon_1=1}{\rightarrow}(\neg f)\\
r_2: b\overset{w_2=1-\epsilon_2=1}\rightarrow f \\
r_3: p\overset{w_3=1-\epsilon_3=1}\rightarrow b
\end{cases}
\end{equation*}

From rules $r_1$ and $r_2$ and because of property \eqref{eq:MI}, one concludes that
$$p\wedge b \rightarrow (f \wedge \neg f)$$
\noindent
and using the non-contradiction law \eqref{eq:NCL} with the Modus Tollens \eqref{eq:MT}, one finally gets
$$\neg (f \wedge \neg f)  \rightarrow \neg(p\wedge b)$$
\noindent which proves that $p\wedge b$ is always false whatever the rule $r_3$ is. Interpreted in terms of the probability theory, the event $T=p\cap b$ corresponds actually and truly to the impossible event $\emptyset$ since $T=f$ and $T=\bar{f}$ are exclusive and exhaustive events. Under such conditions, the analysis proves the non-existence of the penguin-bird Tweety.\\

If one adopts the notations\footnote{Because probabilities are related to sets, we use here the common set-complement notation $\bar{f}$ instead of the logical negation notation $\neg f$, $\cap$ for $\wedge$ and $\cup$ for $\vee$ if necessary.} of the probability theory,  trying to derive $P(T=f | T=p\cap b)$ and $P(T=\bar{f}| T=p\cap b)$ with the Bayesian reasoning is just impossible because from one of the axioms of the probability theory, one must have $P(\emptyset)=0$ and from the  conditioning rule, one would get expressly for this problem the indeterminate expressions:
\begin{align*}
P(T=f| T=p\cap b)&=P(T=f | T=\emptyset)\\
P(T=f| T=p\cap b)           & =\frac{P(T=f \cap \emptyset)}{P(T=\emptyset)}\\
P(T=f| T=p\cap b)           & =\frac{P(T=\emptyset)}{P(T=\emptyset)}\\
P(T=f| T=p\cap b)           & =\frac{0}{0}\qquad\qquad \text{(indeterminate)}
 \end{align*}
\noindent
and similarly      
\begin{align*}
P(T=\bar{f}| T=p\cap b)&=P(T=\bar{f} | T=\emptyset)\\
P(T=\bar{f}| T=p\cap b)           & =\frac{P(T=\bar{f} \cap \emptyset)}{P(T=\emptyset)}\\
P(T=\bar{f}| T=p\cap b)           & =\frac{P(T=\emptyset)}{P(T=\emptyset)}\\
P(T=\bar{f}| T=p\cap b)           & =\frac{0}{0}\qquad\qquad \text{(indeterminate)}
 \end{align*}

\subsubsection{Analysis of the problem when   $0 <\epsilon_1,\epsilon_2,\epsilon_3<1$}

Let's examine now the general case when one allows some little doubt on the rules characterized by taking $\epsilon_1\gtrsim 0$ , $\epsilon_2\gtrsim 0$ and $\epsilon_3\gtrsim 0$ and examine the consequences on the probabilistic model on these rules.\\

First note that, because of the third middle excluded principle and the assumption of the existence of a probabilistic model for a weighted rule, then one should be able to consider simultaneously both "probabilistic/Bayesian" rules
\begin{equation}
\begin{cases}
a \overset{P(b|a)=w}{\rightarrow} b \\
a \overset{P(\bar{b}|a)=1-w}{\rightarrow} \neg b
\end{cases}
 \label{eqW1}
\end{equation}
In terms of classical (objective) probability theory, these weighted rules just indicate that in $100\times w$ percent of cases the logical variable $b$ is true if $a$ is true, or equivalently, that in $100\times w$ percent of cases the random event $b$ occurs when the random event $a$ occurs. When we don't refer to classical probability theory, the weighting factors $w$ and $1-w$ indicate just the level of conviction committed to the validity of the rules. Although very appealing at the first glance, this probabilistic model hides actually a strong drawback/weakness specially when dealing with several rules as shown right below. \\

Let's prove first that from a "probabilized" rule $a \overset{P(b|a)=w}{\rightarrow} b$ one cannot assess rigorously the convictions onto its Modus Tollens. In other words, from \eqref{eqW1} what can we conclude on 
\begin{equation}
\begin{cases}
\neg b \overset{P(\bar{a}|\bar{b})= ?}{\rightarrow} \neg a \\
b  \overset{P(\bar{a}|b)= ?}{\rightarrow} \neg a
\end{cases}
 \label{eqW2}
\end{equation}
\noindent
From the Bayes' rule of conditioning (which must hold if the probabilitic model holds), one can express $P(\bar{a}|\bar{b})$ and $P(\bar{a}|b)$ as follows
\begin{equation*}
\begin{cases}
P(\bar{a}|\bar{b})= 1 - P(a|\bar{b})= 1 - \frac{P(a\cap\bar{b})}{1-P(b)}= 1 -  \frac{P(\bar{b}|a)P(a)}{1-P(b)}\\
P(\bar{a}|b)=1-P(a|b)= 1 - \frac{P(a\cap b)}{P(b)}= 1 - \frac{P(b|a)P(a)}{P(b)}
\end{cases}
 \label{eqW3}
\end{equation*}
\noindent
or equivalently by replacing $P(b|a)$ and $P(\bar{b}|a)$ by their values $w$ and $1-w$, one gets
\begin{equation}
\begin{cases}
P(\bar{a}|\bar{b})= 1 - (1-w) \frac{P(a)}{1-P(b)}\\
P(\bar{a}|b)= 1 - w\frac{P(a)}{P(b)}
\end{cases}
 \label{eqW4}
\end{equation}

\noindent
These relationships show that one cannot fully derive in theory $P(\bar{a}|\bar{b})$ and $P(\bar{a}|b)$ because the prior probabilities $P(a)$ and $P(b)$ are unknown.\\

A simplistic solution, based on the principle of indifference, is then just to assume without solid justification that $P(a)=P(\bar{a})=1/2$ and $P(b)=P(\bar{b})=1/2$. With such assumption, then one gets the following estimates $\hat{P}(\bar{a}|\bar{b})=w$ and $\hat{P}(\bar{a}|b)=1-w$ for $P(\bar{a}|\bar{b})$ and $P(\bar{a}|b)$ respectively and we can go further in the derivations.\\

Now let's go back to our Tweety Penguin Triangle Problem. 
Based on the probabilistic model (assumed to hold), one starts now with both
\begin{equation}
\begin{cases}
r_1: p\overset{P(\bar{f}|p)=1-\epsilon_1}{\rightarrow}\neg f\\
r_2: b\overset{P(f|b)=1-\epsilon_2}\rightarrow f \\
r_3: p\overset{P(b|p)=1-\epsilon_3}\rightarrow b
\end{cases}
\qquad
\begin{cases}
p\overset{P(f|p)=\epsilon_1}{\rightarrow}f\\
b\overset{P(\bar{f}|b)=\epsilon_2}\rightarrow \neg f \\
p\overset{P(\bar{b}|p)=\epsilon_3}\rightarrow \neg b
\end{cases}
\end{equation}
\noindent
Note that taking into account our preliminary analysis and {\it{accepting the principle of indifference}}, one has also the two sets of weighted rules either
\begin{equation}
\begin{cases}
f \overset{\hat{P}(\bar{p}|f)=1-\epsilon_1}{\rightarrow}\neg p\\
\neg f \overset{\hat{P}(\bar{b}|\bar{f})=1-\epsilon_2}\rightarrow \neg b\\
\neg b \overset{\hat{P}(\bar{p}|\bar{b})=1-\epsilon_3}\rightarrow \neg p
\end{cases}
\qquad
\begin{cases}
\neg f \overset{\hat{P}(\bar{p}|\bar{f})=\epsilon_1}{\rightarrow} \neg p\\
f \overset{\hat{P}(\bar{b}|f)=\epsilon_2}\rightarrow \neg b\\
b \overset{\hat{P}(\bar{p}|b)=\epsilon_3}\rightarrow \neg p
\end{cases}
\end{equation}
\noindent
One wants to assess the convictions (assumed to correspond to some conditional probabilities) into the following rules
\begin{equation}p\wedge b \overset{P(f|p\cap b)=?}{\rightarrow} f\end{equation}
\begin{equation}p\wedge b \overset{P(\bar{f}|p\cap b)=?}{\rightarrow} \neg f\end{equation}
The question is to derive rigorously $P(f|p\cap b)$ and $P(\bar{f}|p\cap b)$ from all previous available information. It turns out that the derivation is impossible without unjustified extra assumption on conditional independence. Indeed, $P(f|p\cap b)$  and $P(\bar{f}|p\cap b)$ are given by
\begin{equation}
\begin{cases}
P(f|p\cap b)=\frac{P(f,p, b)}{P(p,b)}=\frac{P(p, b|f)P(f)}{P(b|p)P(p)}\\
\\
P(\bar{f}|p\cap b)=\frac{P(\bar{f},p, b)}{P(p,b)}=\frac{P(p, b|\bar{f})P(\bar{f})}{P(b|p)P(p)}
\end{cases}
\end{equation}

If one assumes as J. Pearl does, that the conditional independence condition also holds, i.e.
$P(p, b|f)=P(p|f)P(b|f)$ and $P(p, b|\bar{f})=P(p|\bar{f})P(b|\bar{f})$, then one gets
\begin{equation*}
\begin{cases}
P(f|p\cap b)=\frac{P(p|f)P(b|f)P(f)}{P(b|p)P(p)}\\
\\
P(\bar{f}|p\cap b)=\frac{P(p|\bar{f})P(b|\bar{f})P(\bar{f})}{P(b|p)P(p)}
\end{cases}
\end{equation*}
\noindent 
By accepting again the principle of indifference, $P(f)=P(\bar{f})=1/2$ and $P(p)=P(\bar{p})=1/2$, one gets the following expressions
\begin{equation}
\begin{cases}
\hat{P}(f|p\cap b)=\frac{P(p|f)P(b|f)}{P(b|p)}\\
\\
\hat{P}(\bar{f}|p\cap b)=\frac{P(p|\bar{f})P(b|\bar{f})}{P(b|p)}
\end{cases}
\label{eq:PP}
\end{equation}
\noindent
Replacing probabilities $P(p|f),P(b|f),P(b|p),P(p|\bar{f})$ and $P(b|\bar{f})$ by their values in the formula \eqref{eq:PP}, one finally gets
\begin{equation}
\begin{cases}
\hat{P}(f|p\cap b)=\frac{\epsilon_1(1-\epsilon_2)}{1-\epsilon_3}\\
\\
\hat{P}(\bar{f}|p\cap b)=\frac{(1-\epsilon_1)\epsilon_2}{1-\epsilon_3}
\end{cases}
\end{equation}

Therefore we see that, even if one accepts the principle of indifference together with the conditional independence assumption, the approximated "probabilities" remain both small and do not correspond to a real measure of probability since the conditional probabilities of exclusive elements $f$ and $\bar{f}$ do not add up to one. When $\epsilon_1$, $\epsilon_2$ and $\epsilon_3$ tends towards 0, one has
$$\hat{P}(f|p\cap b) + \hat{P}(\bar{f}|p\cap b) \approx 0$$
Actually our analysis based on the principle of indifference, the conditional independence assumption and the model proposed by Judea Pearl, proves clearly the impossibility of the Bayesian reasoning to be applied rigorously on such kind of weighted rule-based system, because no probabilistic model exists for describing correctly the problem. This conclusion is actually not surprising taking into account the Lewis' theorem \cite{Lewis_1976} explicated in details in \cite{Goodman_1997} (chapter 11).\\

Let's now explain the reason of the error in the fallacious  reasoning which was looking coherent with the common intuition. The problem arises directly from the fact that penguin class and bird class are defined in this problem only with respect to the "flying" and "not-flying" properties. If one considers only these properties, then none Tweety animal can be categorically classified as a penguin-bird, because penguin-birdness doesn't not hold in reality based on these exclusive and exhaustive properties (if we consider only the information given within the rules $r_1$, $r_2$ and $r_3$). Actually everybody knows that penguins are effectively classified as bird because "birdness" property is not defined with respect to the "flying" or "not-flying" abilities of the animal but by other zoological characteristics $C$ (birds are vertebral oviparous animals with hot blood, a beak, feather and anterior members are wings) and such information must be properly taken into account in the rule-based systems to avoid to fall in the trap of such fallacious  reasoning. The intuition (which seems to justify the fallacious  reasoning conclusion) for TP2 is actually biased because one already knows that penguins (which are truly classified as birds by some other criterions) do not fly in real world and thus we commit a low conviction (which is definitely not a probability measure, but rather a belief) to the fact that a penguin-bird can fly. Thus the Pear'ls analysis proposed in \cite{Pearl_1988} appears to the authors to be unfortunately incomplete and somehow fallacious.

\section{The Dempster-Shafer reasoning}

As pointed out by Judea Pearl in \cite{Pearl_1988},  the Dempster-Shafer reasoning yields, for this problem, a very counter-intuitive result: birdness seems to endow Tweety with extra flying power ! We present here our analysis of this problem based on the Dempster-Shafer reasoning.\\
 
Let's examine in detail  the available prior information summarized by the rule $r_1$: "{\it{Penguins normally don't fly}}" $\Leftrightarrow$ $(p\rightarrow \neg f)$ with the conviction $w_1=1-\epsilon_1$ where $\epsilon_1$ is a small positive number close to zero. This information, in the DST framework, has to be correctly represented in term of a conditional belief  $\text{Bel}_1(\bar{f}|p)=1-\epsilon_1$ rather than directly  the mass $m_1(\bar{f} \cap p)=1-\epsilon_1$.\\
 
 Choosing $\text{Bel}_1(\bar{f}|p)=1-\epsilon_1$ means that there is a high degree of belief that a penguin-animal is also a nonflying-animal (whatever kind of animal we are observing). This representation reflects perfectly our prior knowledge while the erroneous coarse modeling based on the commitment $m_1(\bar{f} \cap p)=1-\epsilon_1$ is unable to distinguish between rule $r_1$ and another (possibly erroneous) rule like $r_1': (\neg f\rightarrow p)$ having same conviction value $w_1$. This correct model allows us to distinguish between $r_1$ and $r_1'$ (even if they have the same numerical level of conviction) by considering the two different conditional beliefs $\text{Bel}_1(\bar{f}|p)=1-\epsilon_1$ and $\text{Bel}_{1'}(p|\bar{f})=1-\epsilon_1$. The coarse/inadequate basic belief assignment modeling (if adopted) in contrary would make no distinction between those two rules $r_1$ and $r_1'$ since one would have to take $m_1(\bar{f} \cap p)=m_{1'}(p \cap \bar{f})$ and therefore cannot serve as the starting model for the analysis\\
 
Similarly, the prior information relative to rules $r_2: (b\rightarrow f)$ and $r_3: (p \rightarrow b)$ with convictions $w_2=1-\epsilon_2$ and  $w_3=1-\epsilon_3$ has to be modeled by the conditional beliefs $\text{Bel}_2(f|b)=1-\epsilon_2$ and $\text{Bel}_3(b|p)=1-\epsilon_3$ respectively.\\
 
The first problem we have to face now is the combination of these three prior information characterized by $\text{Bel}_1(\bar{f}|p)=1-\epsilon_1$, $\text{Bel}_2(f|b)=1-\epsilon_2$ and $\text{Bel}_3(b|p)=1-\epsilon_3$. All the available prior information can be viewed actually as three independent bodies of evidence $\mathcal{B}_1$, $\mathcal{B}_2$ and $\mathcal{B}_3$ providing separately the partial knowledges summarized through the values of $\text{Bel}_1(\bar{f}|p)$, $\text{Bel}_2(f|b)$ and $\text{Bel}_3(b|p)$. To achieve the combination, one needs to define complete basic belief assignments $m_1(.)$, $ m_2(.)$ and $m_3(.)$ compatible with the partial conditional beliefs $\text{Bel}_1(\bar{f}|p)=1-\epsilon_1$, $\text{Bel}_2(f|b)=1-\epsilon_2$ and $\text{Bel}_3(b|p)=1-\epsilon_3$ without introducing extra knowledge. We don't want to introduce  in the derivations some extra-information we don't have in reality. We present in details the justification for the choice of assignment $m_1(.)$. The choice for $m_2(.)$ and $m_3(.)$ will follow similarly.\\

The body of evidence $\mathcal{B}_1$ provides some information only about $\bar{f}$ and $p$ through the value of $\text{Bel}_1(\bar{f}|p)$ and without reference to $b$. Therefore the frame of discernment $\Theta_1$ induced by $\mathcal{B}_1$ and satisfying the Shafer's model (i.e. a finite set of exhaustive and exclusive elements) corresponds to 
$$\Theta_1=\{\theta_1\triangleq \bar{f}\cap\bar{p}, \theta_2\triangleq f\cap\bar{p},\theta_3\triangleq \bar{f}\cap p ,\theta_4\triangleq f\cap p\}$$
\noindent schematically represented by
\begin{equation*}
f=\theta_2\cup\theta_4 \Bigr\{ \overbrace{\underbrace{
\begin{matrix}
\boxed{\theta_4\triangleq f \cap p} &  \boxed{\theta_3\triangleq \bar{f}\cap p}\\
\boxed{\theta_2\triangleq f \cap \bar{p}} &  \boxed{\theta_1\triangleq \bar{f}\cap \bar{p}}
\end{matrix}
}_{\bar{p}=\theta_1\cup\theta_2}}^{p=\theta_3\cup\theta_4} \Bigl \} \bar{f}=\theta_1\cup\theta_3
\end{equation*}
The complete basic assignment $m_1(.)$ we are searching for and defined over the power set $2^{\Theta_1}$ which must be compatible with 
$\text{Bel}_1(\bar{f}|p)$ is actually the result of the Dempster's combination of an unknown (for now) basic belief assignment $m_1'(.)$ with the particular assignment $m_1''(.)$ defined by $m_1''(p\triangleq \theta_3\cup\theta_4)=1$; in other worlds, one has
$$m_1(.)=[m_1'\oplus m_1''](.)$$
From now on, we introduce explicitly  the conditioning term in our notation to avoid confusion and thus we use $m_1(.|p)=m_1(.|\theta_3\cup\theta_4)$ instead $m_1(.)$. From $m_1''(p\triangleq \theta_3\cup\theta_4)=1$ and from any generic unknow basic assignment $m_1'(.)$ defined by its components $m_1'(\emptyset)\triangleq 0$, $m_1'(\theta_1)$,  $m_1'(\theta_2)$,  $m_1'(\theta_3)$,  $m_1'(\theta_4)$, $m_1'(\theta_1\cup\theta_2)$, $m_1'(\theta_1\cup\theta_3)$, $m_1'(\theta_1\cup\theta_4)$, $m_1'(\theta_2\cup\theta_3)$, $m_1'(\theta_2\cup\theta_4)$, $m_1'(\theta_3\cup\theta_4)$, $m_1'(\theta_1\cup\theta_2\cup\theta_3)$, $m_1'(\theta_1\cup\theta_2\cup\theta_4)$, $m_1'(\theta_1\cup\theta_3\cup\theta_4)$, $m_1'(\theta_2\cup\theta_3\cup\theta_4)$, $m_1'(\theta_1\cup\theta_2\cup\theta_3\cup\theta_4)$ and applying Dempter's rule, one gets easily the following expressions for $m_1(.|\theta_3\cup\theta_4)$. All $m_1(.|\theta_3\cup\theta_4)$ masses are zero except theoretically
\begin{align*}
m_1(\theta_3|\theta_3\cup\theta_4) & =  \overbrace{m_1''(\theta_3\cup\theta_4)}^{1} [m_1'(\theta_3) + m_1'(\theta_1\cup\theta_3) \\
 &  \quad+ m_1'(\theta_2\cup\theta_3) \\
 &  \quad+ m_1'(\theta_1\cup\theta_2\cup\theta_3)]/K_1
\end{align*}
\begin{align*}
m_1(\theta_4|\theta_3\cup\theta_4) & =  \overbrace{m_1''(\theta_3\cup\theta_4)}^{1}[m_1'(\theta_4) + m_1'(\theta_1\cup\theta_4) \\
 &  \quad+ m_1'(\theta_2\cup\theta_4) \\
  &  \quad+ m_1'(\theta_1\cup\theta_2\cup\theta_4)]/K_1
\end{align*}
\begin{align*}
m_1(\theta_3\cup\theta_4|\theta_3\cup\theta_4) & =  \overbrace{m_1''(\theta_3\cup\theta_4)}^{1}[m_1'(\theta_3\cup\theta_4) \\
&  \quad+ m_1'(\theta_1\cup\theta_3\cup\theta_4) \\
&  \quad+ m_1'(\theta_2\cup\theta_3\cup\theta_4) \\
&  \quad+ m_1'(\theta_1\cup\theta_2\cup\theta_3\cup\theta_4)]/K_1
\end{align*}
\noindent with $$K_1\triangleq 1 -  \overbrace{m_1''(\theta_3\cup\theta_4)}^{1}[m_1'(\theta_1) + m_1'(\theta_2) + m_1'(\theta_1\cup\theta_2)]$$

To complete the derivation of $m_1(.|\theta_3\cup\theta_4)$, one needs to use the fact that one knows that $\text{Bel}_1(\bar{f}|p)=1-\epsilon_1$ which, by definition \cite{Shafer_1976}, is expressed by
\begin{align*}
\text{Bel}_1(\bar{f}|p) & = \text{Bel}_1(\theta_1\cup\theta_3|\theta_3\cup\theta_4)\\
\text{Bel}_1(\bar{f}|p)& = m_1(\theta_1|\theta_3\cup\theta_4) + m_1(\theta_3|\theta_3\cup\theta_4) \\
&  \quad+ m_1(\theta_1\cup\theta_3|\theta_3\cup\theta_4)\\
\text{Bel}_1(\bar{f}|p)& = 1-\epsilon_1
\end{align*}
But from the generic expression of $m_1(.|\theta_3\cup\theta_4)$, one knows also that $m_1(\theta_1|\theta_3\cup\theta_4)=0$ and $m_1(\theta_1\cup\theta_3|\theta_3\cup\theta_4)=0$. 
Thus the knowledge of $\text{Bel}_1(\bar{f}|p)=1-\epsilon_1$ implies to have
\begin{align*}
m_1(\theta_3|\theta_3\cup\theta_4) & = [m_1'(\theta_3) + m_1'(\theta_1\cup\theta_3) \\
 & \quad + m_1'(\theta_2\cup\theta_3) \\
 &  \quad + m_1'(\theta_1\cup\theta_2\cup\theta_3)]/K_1\\
m_1(\theta_3|\theta_3\cup\theta_4) & = 1-\epsilon_1
 \end{align*}

This is however not sufficient to fully define the values of all components of $m_1(.|\theta_3\cup\theta_4)$ or equivalently of all components of $m_1'(.)$. To complete the derivation without extra unjustified specific information, one needs to apply the minimal commitment principle (MCP) which states that one should never give more support to the truth of a proposition than justified \cite{Klawonn_1992}. According to this principle, we commit a non null value only to the less specific proposition involved into $m_1(\theta_3|\theta_3\cup\theta_4)$ expression. In other words, the MCP allows us to choose legitimately
\begin{align*}
m_1'(\theta_1)&=m_1'(\theta_2)=m_1'(\theta_3)=0\\
m_1'(\theta_1\cup\theta_2)&=m_1'(\theta_1\cup\theta_3)=m_1'(\theta_2\cup\theta_3)=0\\
m_1'(\theta_1\cup\theta_2\cup\theta_3)& \neq 0
 \end{align*}
\noindent
Thus $K_1=1$ and $m_1(\theta_3|\theta_3\cup\theta_4)$ reduces to
\begin{equation*}
m_1(\theta_3|\theta_3\cup\theta_4) =  m_1'(\theta_1\cup\theta_2\cup\theta_3) = 1-\epsilon_1
\end{equation*}

Since the sum of basic belief assignments must be one, one must also have for the remaining (uncommitted for now) masses of $m_1'(.)$ the constraint
\begin{equation*}
\begin{split}
m_1'(\theta_4)+m_1'(\theta_1\cup \theta_4)+ m_1'(\theta_2\cup \theta_4)
 + m_1'(\theta_1\cup \theta_2\cup\theta_4) \\
 + m_1'(\theta_3\cup \theta_4) 
 + m_1'(\theta_1\cup \theta_3\cup\theta_4) +
m_1'(\theta_2\cup \theta_3\cup\theta_4) \\
+ m_1'(\theta_1\cup \theta_2\cup\theta_3\cup\theta_4)=\epsilon_1
\end{split}
\end{equation*}

By applying a second time the MCP, one chooses $m_1'(\theta_1\cup \theta_2\cup\theta_3\cup\theta_4)=\epsilon_1$. \\

Finally, the complete and less specific belief assignment $m_1(.|p)$ compatible with the available prior information $\text{Bel}_1(\bar{f}|p)=1-\epsilon_1$ provided by the source $\mathcal{B}_1$ reduces to
\begin{align}
m_1(\theta_3|\theta_3\cup\theta_4) &=  m_1'(\theta_1\cup\theta_2\cup\theta_3) = 1-\epsilon_1\\
m_1(\theta_3\cup\theta_4|\theta_3\cup\theta_4) &=  m_1'(\theta_1\cup\theta_2\cup\theta_3\cup\theta_4) = \epsilon_1
\end{align}
\noindent
or equivalently
\begin{align}
m_1(\bar{f} \cap p |p) &=  m_1'(\bar{p} \cup \bar{f}) = 1-\epsilon_1\label{eq:(4)}\\
m_1(p|p) &=  m_1'(\bar{p} \cup \bar{f} \cup p \cup f) = \epsilon_1\label{eq:(5)}
\end{align}

It is easy to check, from the mass $m_1(.|p)$, that one gets effectively $\text{Bel}_1(\bar{f}|p)=1-\epsilon_1$. Indeed:
\begin{align*}
\text{Bel}_1(\bar{f}|p) &=\text{Bel}_1(\theta_1\cup\theta_3 |p) \\
\text{Bel}_1(\bar{f}|p)           & = \text{Bel}_1((\bar{f}\cap\bar{p})\cup(\bar{f}\cap p) |p) \\
\text{Bel}_1(\bar{f}|p)           & = \underbrace{m_1(\bar{f}\cap\bar{p}|p)}_{0} + m_1(\bar{f}\cap p |p)\\
          & \qquad  + \underbrace{m_1((\bar{f}\cap\bar{p})\cup(\bar{f}\cap p) |p)}_{0}\\
\text{Bel}_1(\bar{f}|p)           & = m_1(\bar{f}\cap p|p)\\
\text{Bel}_1(\bar{f}|p)           & = 1 -\epsilon_1
\end{align*}
         
In a  similar way, for the source $\mathcal{B}_2$ with $\Theta_2$ defined as
$$\Theta_2=\{\theta_1\triangleq f\cap\bar{b}, \theta_2\triangleq \bar{b}\cap\bar{f},\theta_3\triangleq f\cap b ,\theta_4\triangleq \bar{f}\cap b\}$$
\noindent schematically represented by
\begin{equation*}
\bar{f}=\theta_2\cup\theta_4 \Bigr\{ \overbrace{\underbrace{
\begin{matrix}
\boxed{\theta_4\triangleq \bar{f} \cap b} &  \boxed{\theta_3\triangleq f\cap b}\\
\boxed{\theta_2\triangleq \bar{f}\cap \bar{b}} &  \boxed{\theta_1\triangleq f\cap \bar{b}}
\end{matrix}
}_{\bar{b}=\theta_1\cup\theta_2}}^{b=\theta_3\cup\theta_4} \Bigl \} f=\theta_1\cup\theta_3
\end{equation*}

\noindent
one looks for $m_2(.|b)=[m_2'\oplus m_2''](.)$ with $m_2''(b)=m_2''(\theta_3\cup\theta_4)=1$. From the MCP, the condition $\text{Bel}_2(f|b)=1-\epsilon_2$ and with simple algebraic manipulations, one finally gets
\begin{align}
m_2(\theta_3|\theta_3\cup\theta_4) &=  m_2'(\theta_1\cup\theta_2\cup\theta_3) = 1-\epsilon_2\\
m_2(\theta_3\cup\theta_4|\theta_3\cup\theta_4) &=  m_2'(\theta_1\cup\theta_2\cup\theta_3\cup\theta_4) = \epsilon_2
\end{align}
\noindent
or equivalently
\begin{align}
m_2(f \cap b |b) &=  m_2'(\bar{b} \cup f) = 1-\epsilon_2\label{eq:(8)}\\
m_2(b|b) &=  m_2'(\bar{b} \cup \bar{f} \cup b \cup f) = \epsilon_2\label{eq:(9)}
\end{align}

In a  similar way, for the source $\mathcal{B}_3$ with $\Theta_3$ defined as
$$\Theta_3=\{\theta_1\triangleq b\cap\bar{p}, \theta_2\triangleq \bar{b}\cap\bar{p},\theta_3\triangleq p\cap b ,\theta_4\triangleq \bar{b}\cap p\}$$
\noindent schematically represented by
\begin{equation*}
\bar{b}=\theta_2\cup\theta_4 \Bigr\{ \overbrace{\underbrace{
\begin{matrix}
\boxed{\theta_4\triangleq \bar{b} \cap p} &  \boxed{\theta_3\triangleq b\cap p}\\
\boxed{\theta_2\triangleq \bar{b}\cap \bar{p}} &  \boxed{\theta_1\triangleq b\cap \bar{p}}
\end{matrix}
}_{\bar{p}=\theta_1\cup\theta_2}}^{p=\theta_3\cup\theta_4} \Bigl \} b=\theta_1\cup\theta_3
\end{equation*}

\noindent
one looks for $m_3(.|p)=[m_3'\oplus m_3''](.)$ with $m_3''(p)=m_3''(\theta_3\cup\theta_4)=1$. From the MCP, the condition $\text{Bel}_3(b|p)=1-\epsilon_3$ and with simple algebraic manipulations, one finally gets
\begin{align}
m_3(\theta_3|\theta_3\cup\theta_4) &=  m_3'(\theta_1\cup\theta_2\cup\theta_3) = 1-\epsilon_3\\
m_3(\theta_3\cup\theta_4|\theta_3\cup\theta_4) &=  m_3'(\theta_1\cup\theta_2\cup\theta_3\cup\theta_4) = \epsilon_3
\end{align}
\noindent
or equivalently
\begin{align}
m_3(b \cap p |p) &=  m_3'(\bar{p} \cup b) = 1-\epsilon_3\\
m_3(p|p) &=  m_3'(\bar{b} \cup \bar{p} \cup b \cup p) = \epsilon_3
\end{align}

Since all the complete prior basic belief assignments  are available, one can combine them with the Dempster's rule to summarize all our prior knowledge drawn from our simple rule-based expert system characterized by rules $R=\{r_1,r_2,r_3\}$ and convictions/confidences $W=\{w_1,w_2,w_3\}$ in these rules.\\

The fusion operation requires to primilarily choose the following frame of discernment $\Theta$ (satisfying the Shafer's model) given by 
$$\Theta=\{\theta_1,\theta_2,\theta_3,\theta_4,\theta_5,\theta_6,\theta_7,\theta_8\}$$
\noindent
where
\begin{align*}
&\theta_1 \triangleq f \cap b \cap p \qquad & \theta_5 \triangleq \bar{f }\cap b \cap p\\
&\theta_2 \triangleq f \cap b \cap \bar{p} \qquad &\theta_6 \triangleq \bar{f }\cap b \cap \bar{p}\\
&\theta_3 \triangleq f \cap \bar{b} \cap p \qquad &\theta_7 \triangleq \bar{f }\cap \bar{b} \cap p\\
&\theta_4 \triangleq f \cap \bar{b} \cap \bar{p} \qquad &\theta_8 \triangleq \bar{f }\cap \bar{b} \cap \bar{p}
\end{align*}

The fusion of masses $m_1(.)$ given by eqs. \eqref{eq:(4)}-\eqref{eq:(5)} with $m_2(.)$ given by eqs. \eqref{eq:(8)}-\eqref{eq:(9)} using the Demspter's rule of combination \cite{Shafer_1976} yields $m_{12}(.)=[m_1\oplus m_2](.)$ with the following non null components
\begin{align*}
m_{12}( f \cap b \cap p)&= \epsilon_1(1-\epsilon_2)/K_{12}\\
m_{12}( \bar{f} \cap b \cap p)&= \epsilon_2(1-\epsilon_1)/K_{12}\\
m_{12}(b \cap p)&= \epsilon_1\epsilon_2/K_{12}
\end{align*}
\noindent
with $K_{12}\triangleq 1 - (1-\epsilon_1)(1-\epsilon_2)=\epsilon_1 + \epsilon_2 - \epsilon_1\epsilon_2$.\\

The fusion of all prior knowledge by the Dempster's rule $m_{123}(.)=[m_1\oplus m_2\oplus m_3](.)=[m_{12}\oplus m_3](.)$ yields the final result :
\begin{align*}
m_{123}( f \cap b \cap p)=m_{123}(\theta_1) &= \epsilon_1(1-\epsilon_2)/K_{123}\\
m_{123}( \bar{f} \cap b \cap p)=m_{123}(\theta_5)&= \epsilon_2(1-\epsilon_1)/K_{123}\\
m_{123}(b \cap p)=m_{123}(\theta_1\cup \theta_5)&= \epsilon_1\epsilon_2/K_{123}
\end{align*}
\noindent
with $K_{123}=K_{12}\triangleq 1 - (1-\epsilon_1)(1-\epsilon_2)=\epsilon_1 + \epsilon_2 - \epsilon_1\epsilon_2$.\\

\noindent 
which defines actually and precisely the conditional belief assignment $m_{123}( .|p\cap b)$. It turns out that the fusion with the last basic belief assignment $m_3(.)$ brings no change with respect to previous fusion result $m_{12}(.)$ in this particular problem. \\

Since we are actually interested to assess the belief that our observed particular penguin-animal named Tweety (denoted as $T=(p\cap b)$) can fly, we need to combine all our prior knowledge $m_{123}(.)$ drawn from our rule-based system with the belief assignment $m_o(T=(p\cap b))=1$ characterizing the observation about Tweety. Applying again the Demspter's rule, one finally gets the resulting conditional basic belief function $m_{o123}=[m_o\oplus m_{123}](.)$ defined by
\begin{align*}
m_{o123}( T=(f \cap b \cap p)|T=(p\cap b))&= \epsilon_1(1-\epsilon_2)/K_{12}\\
m_{o123}( T=(\bar{f} \cap b \cap p)|T=(p\cap b))&= \epsilon_2(1-\epsilon_1)/K_{12}\\
m_{o123}(T=(b \cap p)|T=(p\cap b))&= \epsilon_1\epsilon_2/K_{12}
\end{align*}
\noindent
From the Dempster-Shafer reasoning, the belief and plausibity that Tweety can fly are given by \cite{Shafer_1976}
%
\begin{multline*}
\text{Bel}(T=f | T=(p\cap b)) = \\
\sum_{x\in 2^\Theta, 
x\subseteq f } m_{o123}( T=x|T=(p\cap b))
\end{multline*}
\begin{multline*}
\text{Pl}(T=f | T=(p\cap b)) = \\
\sum_{x\in 2^\Theta, 
x\cap f \neq \emptyset} m_{o123}( T=x|T=(p\cap b))
\end{multline*}
\noindent
Because $f=[(f \cap b \cap p)\cup(f \cap b \cap \bar{p})\cup(f \cap \bar{b} \cap p)\cup(f \cap \bar{b} \cap \bar{p})]$ and 
the specific values of the masses defining $m_{o123}(.)$, one has
\begin{multline*}
\text{Bel}(T=f | T=(p\cap b)) = \\
m_{o123}( T=(f \cap b \cap p)|T=(p\cap b))
\end{multline*}
\begin{multline*}
\text{Pl}(T=f | T=(p\cap b)) = \\
m_{o123}( T=(f \cap b \cap p)|T=(p\cap b))\\
+m_{o123}( T=(b \cap p)|T=(p\cap b))
\end{multline*}
\noindent and finally
\begin{equation}
\text{Bel}(T=f | T=(p\cap b)) =  \frac{\epsilon_1(1-\epsilon_2)}{K_{12}}
\end{equation}
\begin{equation}
\text{Pl}(T=f | T=(p\cap b)) =  \frac{\epsilon_1(1-\epsilon_2)}{K_{12}}+\frac{\epsilon_1\epsilon_2}{K_{12}}=\frac{\epsilon_1}{K_{12}}
\end{equation}
In a similar way, one will get for the belief and the plausibility that Tweety cannot fly
\begin{equation}
\text{Bel}(T=\bar{f} | T=(p\cap b)) =  \frac{\epsilon_2(1-\epsilon_1)}{K_{12}}
\end{equation}
\begin{equation}
\text{Pl}(T=\bar{f} | T=(p\cap b)) =  \frac{\epsilon_2(1-\epsilon_1)}{K_{12}}+\frac{\epsilon_1\epsilon_2}{K_{12}}=\frac{\epsilon_2}{K_{12}}
\end{equation}

Using the first order approximation when $\epsilon_1$ and $\epsilon_2$ are very small positive numbers, one gets finally

\begin{equation*}
\text{Bel}(T=f | T=(p\cap b))=\text{Pl}(T=f | T=(p\cap b)) \approx  \frac{\epsilon_1}{\epsilon_1+\epsilon_2}
\end{equation*}
In a similar way, one will get for the belief that Tweety cannot fly
\begin{equation*}
\text{Bel}(T=\bar{f} | T=(p\cap b))=\text{Pl}(T=\bar{f} | T=(p\cap b)) \approx  \frac{\epsilon_2}{\epsilon_1+\epsilon_2}
\end{equation*}

\noindent 
This result coincides with the Judea Pearl's result but a different analysis and detailed presentation has been done here.  
It turns out that this simple and complete analysis corresponds actually to the ballooning extension and the generalized Bayesian theorem proposed by Smets in \cite{Smets_1978,Smets_1993} and discussed by Shafer in \cite{Shafer_1982} although it was carried out independently of Smets' works.
As pointed out by Judea Pearl, this result based on DST and the Dempster's rule of combination looks very paradoxical/counter-intuitive since it means that if nonflying birds are very rare, i.e. $\epsilon_2 \approx 0$, then penguin-birds like our observed penguin-bird Tweety, have a very big chance of flying. As stated by Judea Pearl in  \cite{Pearl_1988} pages 448-449: {\it{''The clash with intuition revolves not around the exact numerical value of $\text{Bel}(f)$ but rather around the unacceptable phenomenon that rule $r_3$, stating that penguins are a subclass of birds, plays no role in the analysis. Knowing that Tweety is both a penguin and a bird renders $\text{Bel}(T=f | T=(p\cap b))$ solely a function of $m_1(.)$ and $m_2(.)$, regardless of how penguins and birds are related. This stands contrary to common discourse, where people expect class properties to be overridden by properties of more specific subclasses. While in classical logic the three rules in our example would yield an unforgivable contradiction, the uncertainties attached to these rules, together with Dempster's normalization, now render them manageable. However, they are managed in the wrong way whenever we interpret if-then rules as randomized logical formulas of the material-implication type, instead of statements of conditional probabilities''}}. Keep in mind that this Pearl's statement is however given to show the semantic clash between the Dempster-Shafer reasoning vs. the fallacious Bayesian reasoning  to support the Bayesian reasoning approach.

\section{The Dezert-Smarandache reasoning}

Before going further in our analysis, some clarification is necessary to explain to the reader the fundamental difference between the foundations of DSmT vs. DST. The DSmT can be easily viewed as a general flexible Bottom-Up approach for managing uncertainty and conflicts in fusion problems. It arises from the fact that the conflict between sources of evidence can come not only from the reliability of sources themselve (which can be handled quite easily by classical discounting methods) but also from a different interpretation of elements of the frame just because the sources or evidence have only a limited knowlege and provide their beliefs only with respect to their knowledge based usually on their own (local) experience, not to mention the fact that elements of the frame of the problem can truly be not refinable at all in some cases involving vague concepts like smallness/tallness, pleasure/pain, etc because of the continuous path from one to the other, etc. Based on this matter of fact, the DSmT proposes a new mathematical framework which starts at the bottom level (solid ground level) from the free DSm model and the notion of hyper-power set (Dedekind's lattice), then provides a general rule of combination to work with the free DSm model. Then it includes the possibility to take into account any kind of integrity constraints into the free DSm model if necessary through the hybrid DSm rule of combination. The taking into account for an integrity constraint consists just in forcing some elements of the Dedekind's lattice to be empty, just because they truly are for some given problems.\\

 The introduction of an integrity constraint is like "pushing an elevator button" for going a bit up in the process of managing uncertainty and conflicts. If one needs to go higher, then one can take into account several integrity constraints as well in the framework of DSmT. If we finally wants to take into account all possible exclusivity constraints if we know that all elements of the frame of the given problem under consideration are truly exclusive, then we go directly to the Top level (the Shafer's model which serves as foundation for the DST). \\
 
 DSmT however can handle not only exclusivity constraints, but also existential constraints or mixed constraints as well which is helpful for some dynamic fusion problems. 
It is also important to emphaze that the hybrid DSm rule of combination is definitely not equivalent to the Dempster's rule of combination (and its alternatives based on the Top level) because one can stop and work at any level in the process of managing uncertainty and conflicts, depending on the nature of the problem. The hybrid DSm rule and Dempster's rule do not provide same results even if working with the Shafer's model as it will be proved in the sequel. The approach proposed by the DSmT to attack the fusion problem is totally new both by its foundations and the solution provided.\\

The DSmT has been originally (ground-level) developed for the fusion of uncertain and paradoxical (highly conflicting) sources of information (bodies of evidences) based on the free DSm model $\mathcal{M}^f(\Theta)$ which assumes that none of elements of the frame $\Theta$ are exclusive. This model is opposite to the Shafer's model. Let consider a free DSm model $\mathcal{M}^f(\Theta)$ with $\Theta=\{\theta_{1},\ldots,\theta_{n}\}$, the DSmT starts with the notion of hyper-power set $D^\Theta$ defined as the set of all composite propositions built from elements of $\Theta$ with $\cup$ and $\cap$ ($\Theta$ generates $D^\Theta$ under operators $\cup$ and $\cap$)
operators such that  \cite{Dezert_2003f}
\begin{enumerate}
\item $\emptyset, \theta_1,\ldots, \theta_n \in D^\Theta$.
\item  If $A,B \in D^\Theta$, then $A\cap B\in D^\Theta$ and $A\cup B\in D^\Theta$.
\item No other elements belong to $D^\Theta$, except those obtained by using rules 1 or 2.
\end{enumerate}
The cardinality of hyper-power set, $d(n)\triangleq\vert D^\Theta\vert$ for $n\geq1$, follows the sequence of Dedekind's numbers 1, 2, 5, 19, 167, 7580, 7828353, ... More details about the generation and partial ordering of elements of hyper-power set can be found in \cite{Dezert_Smarandache_2003,Dezert_2003f,Dezert_2004Book}. From this model, authors have proposed a new simple associative and commutative rule of combination (the DSm classic rule) and then extended this rule to deal with any kind of hybrid models, i.e. sets $\Theta$ for which some propositions/elements of $D^\Theta$ are known or forced to be empty depending on the nature and the dynamicity of the fusion problem under consideration. In this framework, the Shafer's model appears only as a special hybrid model (the most constrained one, if we don't introduce existential constraints). The hybrid DSm fusion rule covers a wide class of fusion applications but is restricted to fusion of precise uncertain and paradoxical information only \cite{Dezert_2004Book}. We have recently extended this rule with new set operators for the fusion of imprecise, uncertain and paradoxical information - see \cite{Dezert_2004Book} for details.\\

We analyze here the Tweety penguin triangle problem with the DSmT.
The prior knowledge characterized by the rules $R=\{r_1,r_2,r_3\}$ and convictions $W=\{w_1,w_2,w_3\}$ is modeled as three independent sources of evidence defined on separate minimal and potentially paradoxical (i.e internal conflicting) frames $\Theta_1\triangleq\{ p,\bar{f}\}$, $\Theta_2\triangleq\{b,f\}$ and $\Theta_3\triangleq\{p,b\}$ since the rule $r_1$ doesn't refer to the existence of $b$, the rule $r_2$ doesn't refer to the existence of $p$ and the rule $r_3$ doesn't refer to the existence of $f$ or $\bar{f}$.  Let's note that the DSmT doesn't require the refinement of frames as with DST (see previous section). We follow the same analysis as in previous section but  now based on our DSm reasoning and the DSm rule of combination.\\

The first source $\mathcal{B}_1$ relative to $r_1$  with confidence $w_1=1-\epsilon_1$ provides us the conditional belief $\text{Bel}_1(\bar{f}|p)$ which is now defined from a paradoxical basic belief assignment $m_1(.)$ resulting from the DSm combination of $m_1''(p)=1$ with $m_1'(.)$ defined on the hyper-power set $D^{\Theta_1}=\{\emptyset,p,\bar{f},p\cap\bar{f}, p\cup\bar{f}\}$. The choice for $m_1'(.)$ results directly from the derivation of the DSm rule and the application of the MCP. Indeed, the non null components of $m_1(.)$ are given by (we introduce explicitly the conditioning term in notation for convenience):
\begin{align*}
m_1(p|p)&=\overbrace{m_1''(p)}^{1}m_1'(p) +\overbrace{m_1''(p)}^{1}m_1'(p\cup \bar{f}) \\
m_1(p\cap \bar{f}|p) &=\overbrace{m_1''(p)}^{1}m_1'(\bar{f}) + \overbrace{m_1''(p)}^{1}m_1'(p\cap\bar{f})\\
\end{align*}

The information $\text{Bel}_1(\bar{f}|p)=1-\epsilon_1$ implies 

$$\text{Bel}_1(\bar{f}|p)= m_1(\bar{f}|p) + m_1(p\cap \bar{f}|p)=1-\epsilon_1$$

Since $m_1(p|p)+ m_1(p\cap \bar{f}|p) =1$, one has necessarily $m_1(\bar{f}|p)=0$ and thus from previous equation $m_1(\bar{f}\cap p|p)=1-\epsilon_1$, which implies both

\begin{align*}
m_1(p|p) & = \epsilon_1\\
m_1(p\cap \bar{f}|p) & =\overbrace{m_1''(p)}^{1}m_1'(\bar{f}) + \overbrace{m_1''(p)}^{1}m_1'(p\cap\bar{f})\\
 & =m_1'(\bar{f}) +m_1'(p\cap\bar{f})\\
 & = 1-\epsilon_1
\end{align*}
 
 Applying the MCP, it results that one must choose
 $$m_1'(\bar{f})=1-\epsilon_1 \quad \text{and} \quad m_1'(p\cap\bar{f})=0$$

The sum of remaining masses of $m_1'(.)$ must be then equal to $\epsilon_1$, i.e.

$$m_1'(p) + m_1'(p\cup\bar{f}) = \epsilon_1$$

Applying again the MCP on this last constraint, one gets naturally
$$m_1'(p)=0 \quad \text{and} \quad  m_1'(p\cup\bar{f}) = \epsilon_1$$

Finally the belief assignment $m_1(.|p)$ relative to the source $\mathcal{B}_1$ and compatible with the constraint $\text{Bel}_1(\bar{f}|p)=1-\epsilon_1$, holds the same numerical values as within the DST analysis (see eqs. \eqref{eq:(4)}-\eqref{eq:(5)})  and is given by
\begin{align*}
m_1(p\cap \bar{f}|p) & =1-\epsilon_1\\
m_1(p|p) & = \epsilon_1
\end{align*}
\noindent 
but results here from the DSm combination of the two following assignments (i.e. $m_1(.)=[m_1'\oplus m_1''](.)=[m_1''\oplus m_1'](.)$)
\begin{equation}
\begin{cases}
m_1'(\bar{f})=1-\epsilon_1 \quad \text{and} \quad m_1'(p\cup\bar{f}) = \epsilon_1\\
m_1''(p)=1
\end{cases}
\end{equation}

In a similarly manner and working on $\Theta_2=\{b,f\}$ for source  $\mathcal{B}_2$ with the condition $\text{Bel}_2(f|b)=1-\epsilon_2$, the mass $m_2(.|b)$ results from the internal DSm combination of the two following assignments
\begin{equation}
\begin{cases}
m_2'(f)=1-\epsilon_2 \quad \text{and} \quad m_2'(b\cup f) = \epsilon_2\\
m_2''(b)=1
\end{cases}
\end{equation}

Similarly and working on $\Theta_3=\{p,b\}$ for source  $\mathcal{B}_3$ with the condition $\text{Bel}_3(b|p)=1-\epsilon_3$, the mass $m_3(.|p)$ results from the internal DSm combination of the two following assignments
\begin{equation}
\begin{cases}
m_3'(b)=1-\epsilon_3 \quad \text{and} \quad m_3'(b\cup p) = \epsilon_3\\
m_3''(p)=1
\end{cases}
\end{equation}

It can be easily verified that these (less specific) basic belief assignments generates the conditions
$\text{Bel}_1(\bar{f}|p)=1-\epsilon_1$, $\text{Bel}_2(f|b)=1-\epsilon_2$ and $\text{Bel}_3(b|p)=1-\epsilon_3$.\\

Now let's examine the result of the fusion of all these masses based on DSmT, i.e by applying the DSm rule of combination of the following basic belief assignments 
$$m_1(p\cap\bar{f}|p)=1-\epsilon_1 \quad\text{and}\quad m_1(p|p)=\epsilon_1$$
$$m_2(b\cap f|b)=1-\epsilon_2 \quad\text{and}\quad m_2(b|b)=\epsilon_2$$
$$m_3(p\cap b|p)=1-\epsilon_3 \quad\text{and}\quad m_3(p|p)=\epsilon_3$$

Note that these basic belief assignments turn to be identical to those drawn from DST framework analysis done in previous section for this specific problem because of  integrity constraint $f\cap\bar{f}=\emptyset$ and the MCP, but result actually from a slightly different and simpler analysis here drawn from DSmT. So we attack the TP2 with the same information as with the analysis based on DST, but we will show that a coherent conclusion can be drawn with DSm reasoning.\\

Let's emphasize now that one has to deal here with the hypotheses/elements $p$, $b$, $f$ and $\bar{f}$ and thus our global frame is given by $\Theta=\{b,p,f,\bar{f}\}$. Note that $\Theta$ doesn't satisfy the Shafer's model since the elements of $\Theta$ are not all exclusive. This is a major difference between the foundations of DSmT with respect to the foundations of DST.  But because only $f$ and $\bar{f}$ are truly exclusive, i.e. $\bar{f}\cap f=\emptyset$, we face a simple hybrid DSm model $\mathcal{M}$ and thus the hybrid DSm fusion must apply rather than the classic DSm rule. We recall briefly here (a complete derivation, justification and examples can be found in \cite{Dezert_2004Book}) the hybrid DSm rule of combination associated to a given hybrid DSm model for $k\geq 2$ independent sources of information is defined for all $A\in D^\Theta$ as:
\begin{equation}
m_{\mathcal{M}(\Theta)}(A)\triangleq 
\phi(A)\Bigl[ S_1(A) + S_2(A) + S_3(A)\Bigr]
 \label{eq:DSmHkBis}
\end{equation}
\noindent
where $\phi(A)$ is the {\it{characteristic emptiness function}} of the set $A$, i.e. $\phi(A)= 1$ if  $A\notin \boldsymbol{\emptyset}$ ($\boldsymbol{\emptyset}\triangleq \{\emptyset,\boldsymbol{\emptyset}_{\mathcal{M}}\}$ being the set of all relatively and absolutely empty elements) and $\phi(A)= 0$ otherwise, and
\begin{equation}
S_1(A)\triangleq \sum_{\substack{X_1,X_2,\ldots,X_k\in D^\Theta \\ (X_1\cap X_2\cap\ldots\cap X_k)=A}} \prod_{i=1}^{k} m_i(X_i)
\end{equation}
\begin{equation}
S_2(A)\triangleq \sum_{\substack{X_1,X_2,\ldots,X_k\in\boldsymbol{\emptyset} \\ [\mathcal{U}=A]\vee [(\mathcal{U}\in\boldsymbol{\emptyset}) \wedge (A=I_t)]}} \prod_{i=1}^{k} m_i(X_i)\end{equation}
\begin{equation}
S_3(A)\triangleq\sum_{\substack{X_1,X_2,\ldots,X_k\in D^\Theta \\ (X_1\cup X_2\cup\ldots\cup X_k)=A \\ (X_1\cap X_2\cap \ldots\cap X_k)\in\boldsymbol{\emptyset}}}  \prod_{i=1}^{k} m_i(X_i)
\end{equation}
with $\mathcal{U}\triangleq u(X_1)\cup u(X_2)\cup \ldots \cup u(X_k)$ where $u(X)$ is the union of all singletons $\theta_i$ that compose $X$ and $I_t\triangleq\theta_1\cup \theta_2\cup \ldots\cup \theta_n$ is the total ignorance defined on the frame $\Theta=\{\theta_1,\ldots,\theta_n\}$.  For example, if $X$ is a singleton then $u(X)=X$; if  $X=\theta_1\cap \theta_2$ or $X=\theta_1\cup \theta_2$ then $u(X)=\theta_1\cup \theta_2$;
if $X=(\theta_1\cap \theta_2)\cup \theta_3$ then $u(X)=\theta_1\cup \theta_2\cup\theta_3$; by convention $u(\emptyset)\triangleq\emptyset$.\\

The first sum $S_1(A)$ entering in the previous formula corresponds to mass $m_{\mathcal{M}^f(\Theta)}(A)$ obtained by the classic DSm rule of combination based on the free DSm model $\mathcal{M}^f$ (i.e. on the free lattice $D^\Theta$). The second sum $S_2(A)$ entering in the formula of the hybrid DSm rule of combination \eqref{eq:DSmHkBis} represents the mass of all relatively and absolutely empty sets which is transferred to the total or relative ignorances. The third sum $S_3(A)$ entering in the formula of the hybrid DSm rule of combination \eqref{eq:DSmHkBis} transfers the sum of relatively empty sets to the non-empty sets in the same way as it was calculated following the DSm classic rule.\\

To apply the DSm hybrid fusion rule formula \eqref{eq:DSmHkBis}, it is important to note that 
$(p\cap \bar{f})\cap (b\cap f)\cap p \equiv p \cap b \cap f \cap \bar{f} = \emptyset$ because $f\cap \bar{f}=\emptyset$, thus the mass $(1-\epsilon_1)(1-\epsilon_2)\epsilon_3$ is transferred to the hybrid proposition $H_1\triangleq (p\cap \bar{f})\cup (b\cap f)\cup p \equiv (b \cap f )\cup p$; similarly $(p\cap \bar{f})\cap (b\cap f)\cap (p\cap b) \equiv p \cap b \cap f \cap \bar{f} = \emptyset$ because $f\cap \bar{f}=\emptyset$ and therefore its associated mass $(1-\epsilon_1)(1-\epsilon_2)(1-\epsilon_3)$ is transferred to the hybrid proposition $H_2\triangleq (p\cap \bar{f})\cup (b\cap f)\cup (p\cap b)$. No other mass transfer is necessary for this Tweety Penguin Triangle Problem and thus we finally get from DSm hybrid fusion formula \eqref{eq:DSmHkBis} the following result for $m_{123}(.|p\cap b)=[m_1\oplus m_2 \oplus m_3](.)$ (where $\oplus$ symbol corresponds here to the DSm fusion operator):
\begin{align*}
m_{123}(H_1|p\cap b) & = (1-\epsilon_1)(1-\epsilon_2)\epsilon_3\\
m_{123}(H_2|p\cap b) & = (1-\epsilon_1)(1-\epsilon_2)(1-\epsilon_3)\\
m_{123}(p\cap b \cap \bar{f}|p\cap b) & = (1-\epsilon_1)\epsilon_2\epsilon_3 +   (1-\epsilon_1)\epsilon_2(1-\epsilon_3)\\
m_{123}(p\cap b \cap f|p\cap b) & = \epsilon_1(1-\epsilon_2)\epsilon_3 +  \epsilon_1(1-\epsilon_2)(1-\epsilon_3)\\
m_{123}(p\cap b|p\cap b) & = \epsilon_1\epsilon_2\epsilon_3 + \epsilon_1\epsilon_2(1-\epsilon_3)
\end{align*}
\noindent
with 
$$\begin{cases}
H_1\triangleq (b \cap f )\cup p\\
H_2\triangleq (p\cap \bar{f})\cup (b\cap f)\cup (p\cap b)
\end{cases}
$$

It can be easily checked that these masses sum up to 1. After elementary algebraic simplifications, one finally gets for the DSm fusion of all available prior information and reintroducing explicitly the conditioning term
\begin{align*}
m_{123}(H_1|p\cap b) & = (1-\epsilon_1)(1-\epsilon_2)\epsilon_3\\
m_{123}(H_2|p\cap b) & = (1-\epsilon_1)(1-\epsilon_2)(1-\epsilon_3)\\
m_{123}(p\cap b \cap \bar{f}|p\cap b) & = (1-\epsilon_1)\epsilon_2\\
m_{123}(p\cap b \cap f|p\cap b) & = \epsilon_1(1-\epsilon_2)\\
m_{123}(p\cap b|p\cap b) & = \epsilon_1\epsilon_2
\end{align*}

We can check all these masses add up to 1 and that this result is fully coherent with the {\it{rational intuition}} specially when $\epsilon_3=0$, because non null components of $m_{123}(.|p\cap b)$ reduces to
\begin{align*}
m_{123}(H_2|p\cap b) & = (1-\epsilon_1)(1-\epsilon_2)\\
m_{123}(p\cap b \cap \bar{f}|p\cap b) & = (1-\epsilon_1)\epsilon_2\\
m_{123}(p\cap b \cap f|p\cap b) & = \epsilon_1(1-\epsilon_2)\\
m_{123}(p\cap b|p\cap b) & = \epsilon_1\epsilon_2
\end{align*}
which means that from our DSm reasoning there is a strong  uncertainty (due to the conflicting rules of our rule-based system), when $\epsilon_1$ and $\epsilon_2$ remain small positive numbers, that a penguin-bird animal is either a penguin-nonflying animal or a bird-flying animal. The small value $\epsilon_1\epsilon_2$ for $m_{123}(p\cap b|p\cap b)$ expresses adequately the fact that we cannot commit a strong basic belief assignment only to $p\cap b$ knowing $p\cap b$ just because one works on $\Theta=\{p,b,f,\bar{f}\}$ and we cannot consider the property $p\cap b$ solely because the "birdness" or "penguinness" property endow necessary either the flying or non-flying property.\\

Therefore the belief that the particular observed penguin-bird animal Tweety (corresponding to the particular mass $m_o(T=(p\cap b))=1$) can be easily derived from the DSm fusion of all our prior summarized by $m_{123}(.|p\cap b)$ and the available observation summarized by $m_o(.)$ and we get
\begin{align*}
m_{o123}(T=(p\cap b \cap \bar{f})|T=(p\cap b)) & = (1-\epsilon_1)\epsilon_2\\
m_{o123}(T=(p\cap b \cap f)|T=(p\cap b)) & = \epsilon_1(1-\epsilon_2)\\
m_{o123}(T=(p\cap b)|T=(p\cap b)) & = \epsilon_1\epsilon_2\\
m_{o123}(T=H_1|T=(p\cap b)) & = (1-\epsilon_1)(1-\epsilon_2)\epsilon_3\\
m_{o123}(T=H_2|T=(p\cap b)) & = (1-\epsilon_1)(1-\epsilon_2)(1-\epsilon_3)
\end{align*}

\noindent
From the DSm reasoning, the belief that Tweety can fly is then given by
\begin{equation*}
\text{Bel}(T=f | T=(p\cap b)) = 
\sum_{x\in D^\Theta, 
x\subseteq f } m_{o123}( T=x|T=(p\cap b))
\end{equation*}
\noindent
Using all the components of $m_{o123}(.|T=(p\cap b))$, one directly gets
\begin{equation*}
\text{Bel}(T=f | T=(p\cap b)) = 
m_{o123}( T=(f \cap b \cap p)|T=(p\cap b))
\end{equation*}
\noindent and finally
\begin{equation}
\text{Bel}(T=f | T=(p\cap b)) =  \epsilon_1(1-\epsilon_2)
\end{equation}
In a similar way, one will get for the belief that Tweety cannot fly
\begin{equation}
\text{Bel}(T=\bar{f} | T=(p\cap b)) =  \epsilon_2(1-\epsilon_1)
\end{equation}

So now for both cases the beliefs remain very low which is normal and coherent with analysis done in section \ref{SEC:PEARLWEAKNESS}.
Now let's examine the plausibilities of the ability for Tweety to fly or not to fly. These are given by 

\begin{equation*}
\text{Pl}(T=f | T=(p\cap b)) \triangleq 
\sum_{x\in D^\Theta, 
x \cap f \neq \ } m_{o123}( T=x|T=(p\cap b))
\end{equation*}
\begin{equation*}
\text{Pl}(T=\bar{f}| T=(p\cap b)) \triangleq 
\sum_{x\in D^\Theta, 
x \cap \bar{f} \neq \ } m_{o123}( T=x|T=(p\cap b))
\end{equation*}
\noindent
which turn to be after elementary algebraic manipulations
\begin{equation}
\text{Pl}(T=f | T=(p\cap b)) = (1-\epsilon_2)
\end{equation}
\begin{equation}
\text{Pl}(T=\bar{f}| T=(p\cap b)) = (1-\epsilon_1)
\end{equation}

So we conclude, as expected, that we can't decide on the ability for Tweety of flying or of not flying, since one has
$$[\text{Bel}(f |p\cap b),\text{Pl}(f| p\cap b)] =  [\epsilon_1(1-\epsilon_2),(1-\epsilon_2)] \approx [0,1]$$
$$[\text{Bel}(\bar{f} | p\cap b),\text{Pl}(\bar{f}| p\cap b)] =  [\epsilon_2(1-\epsilon_1),(1-\epsilon_1)] \approx [0,1]$$

Note that when setting $\epsilon_1=0$ and  $\epsilon_2=1$ (or $\epsilon_1=1$ and  $\epsilon_2=0$), i.e. one forces the full consistency of the initial rules-based system, one gets coherent result on the certainty of the ability of Tweety to not fly (or to fly respectively).\\

This coherent result (radically different from the one based on Dempster-Shafer reasoning but starting with exactly the same available information) comes from the DSm hybrid fusion rule which  transfers some parts of the mass of empty set $m(\emptyset)=(1-\epsilon_1)(1-\epsilon_2)\epsilon_3 + (1-\epsilon_1)(1-\epsilon_2)(1-\epsilon_3)\approx 1$ onto propositions $H_1$ and $H_2$. It is clear however that the high value of $m(\emptyset)$ in this TP2 indicates a high conflicting fusion problem which proves that the TP2 is a truly almost impossible problem and the fusion result based on DSmT reasoning allows us to conclude on  the true undecidability on the ability for Tweety of flying or of not flying. In other words,  the fusion based on DSmT can be applied adequately on this almost impossible problem and concludes correctly on its undecidability. Another simplistic solution would consist to say naturally that the problem has to be considered as an impossible one just because $m(\emptyset)\ge 0.5$.

\section{ Conclusion}

In this paper we have proposed a deep analysis of the challenging Tweety Penguin Triangle Problem.
The analysis proves that the Bayesian reasoning cannot be mathematically justified to characterize the problem because the probabilistic model doesn't hold, even with the help of acceptance of the principle of indifference and the conditional independence assumption. Any conclusions drawn from such representation of the problem based on a hypothetical probabilistic model are based actually on a fallacious Bayesian reasoning. This is a fundamental result. Then one has shown how the Dempster-Shafer reasoning manages in what we feel is a wrong way the uncertainty and the conflict in this problem. We then proved that the DSmT can deal properly with this problem and provides a well-founded and reasonable conclusion about the undecidability of its solution. 

\section*{Acknowledgments}

Authors are grateful to Dr. Roy Streit, Naval Undersea Warfare Center, Newport, RI, U.S.A., to have them introduced and encouraged to work on this nice exciting and challenging problem during the Fusion 2003 International Conference on Information Fusion, Cairns, Australia, July 2003. Authors want to thank Professor Smets for his discussions, suggestions and for pointing out important references for the improvment of this paper.


\begin{thebibliography}{99}

\bibitem{Benferhat_2000}
Benferhat S., Saffioti A., Smets Ph., \emph{Belief functions and default reasoning}, Artificial Intelligence 122, pp. 1-69, 2000.
\bibitem{Buchanan_1984}
Buchanan B.G., Shortliffe E.H., \emph{Rule Bases Expert Systems - The MYCIN Experiments of the Stanford Heuristic Programming Project}, Addison Wesley, 1984.
\bibitem{Dambreville_2004}
Dambreville F., \emph{Probabilized logics related to DSmT and Bayes inference}, in Advances and Applications of DSmT for Information Fusion (Collected works), Chapter 8, Smarandache F., Dezert J. (Editors), American Research Press, June 2004.
\bibitem{Dezert_Smarandache_2003}
Dezert J., Smarandache F., \emph{Partial ordering of hyper-power sets and matrix representation of belief functions within DSmT}, Proc. of Fusion 2003 Conf., Cairns, Australia, July 8-11, 2003.
\bibitem{Dezert_2003f}
Dezert J., Smarandache F., \emph{On the generation of hyper-power sets for the DSmT}, Proceedings of the 6th International Conference on Information Fusion, Cairns, Australia, July 8-11, 2003.
\bibitem{Dubois_1980}
Dubois D., Prade H., \emph{Fuzzy sets and systems, theory and applications}, Academic Press, 1980.
\bibitem{Dubois_1987}
Dubois D., Prade H., \emph{ThŽorie des possibilitŽs, application ˆ la reprŽsentation des connaissances en informatique}, Masson, (2nd edition), 1987.
\bibitem{Goodman_1997}
Goodman I.R., Mahler R.P.S., Nguyen H.T., \emph{Mathematics of Data Fusion}, Kluwer Acad. Press, 1997.
\bibitem{Klawonn_1992}
Klawonn F., Smets Ph., \emph{The Dynamic of Belief in the Transferable Belief Model and Specialization-Generalization Matrices}, Proc. of the 8th Conf. on Uncertainty in AI. Dubois D., Wellman M.P., DÕAmbrosio B. and Smets Ph. (eds). Morgan Kaufmann Publ., San Mateo, CA,  pp. 130-137, 1992.
\bibitem{Kleene_2002}
Kleene, S.ÊC., \emph{Mathematical Logic}, New York,Dover, 2002. 
\bibitem{Mendelson_1997}
Mendelson, E.,\emph{The Propositional Calculus}, Ch.Ê1 in Introduction to Mathematical Logic, 4th ed. London: Chapman \& Hall, pp.Ê12-44, 1997. 
\bibitem{Pearl_1988}
Pearl J., \emph{Probabilistic reasoning in Intelligent Systems: Networks of Plausible Inference}, Morgan Kaufmann Publishers, San Mateo, CA, 1988.
%
\bibitem{Pearl_1990}
Pearl J., \emph{Reasoning with belief functions and analysis of compatibility}, Int. Journal of Approximate Reasoning, 4, pp. 363-390, 1990.
%
\bibitem{Lewis_1976}
Lewis D., \emph{Probabilities of conditonals and conditional probabilities}, Philsophical Review, 85, 1976.
\bibitem{Poole_1985}
Poole D.L., \emph{On the comparison of theories: preferring the most specific explanation}, Proc. 9th Int. Joint Conf. on Artif. Intelligence, Los Angeles, pp. 465-474, 1985.
\bibitem{Nidditch_1962}
Nidditch, P.ÊH,  \emph{Propositional Calculus}, New York, Free Press of Glencoe, 1962. 
\bibitem{Shafer_1976}
Shafer G., \emph{A Mathematical Theory of Evidence}, Princeton Univ. Press, Princeton, NJ, 1976.
\bibitem{Shafer_1982}
Shafer G., \emph{Belief functions and parametric models}, J. Roy. Statist. Soc., B., 44, pp. 322-352, 1982.
\bibitem{Smarandache_1979}
Smarandache F., \emph{Deducibility Theorems in Mathematics logics},
Ann. Univ. Timisoara, Seria ST., Matematica, Vol. XVII, Fasc. 2, pp. 163-168, 1979 (available in \cite{Smarandache_1995})
\bibitem{Smarandache_1995}
Smarandache F., \emph{Collected Papers, Vol. I}, Editura Tempus Romania SRL, 1995, pp. 232-239,
available at {\scriptsize{\verb+http://www.gallup.unm.edu/~smarandache/CP1.pdf+}}.
\bibitem{Dezert_2004Book}
Smarandache F., Dezert J. (Editors), \emph{Advances and Applications of DSmT for Information Fusion (Collected works)}, American Research Press, June 2004.
\bibitem{Smets_1978}
Smets Ph., \emph{Un mod\`ele math\'ematico-statistique simulant le processus de diagnostic m\'edical}, Ph.D. Thesis, Universit\'e Libre de Bruxelles (available through University Microfilm International, 30-32 Mortimer Street London, W1N 7RA, Thesis 80-70,003), 1978.
\bibitem{Smets_1991}
Smets Ph., Hsia Y.T., \emph{Default reasoning and the transferable belief model}, In M. Henrion, L.N. Kanal and J.F. Lemmer (Editors), Uncertainty in Artificial Intelligence 5, Elsevier Science Publishers, pp. 495-504, 1990.
\bibitem{Smets_1992}
Smets Ph., \emph{Resolving misunderstandings about belief functions: A response to the many criticisms raised by Judea Pearl}, Int. J. Approximate Reasoning, 6, pp. 321-344, 1992.
\bibitem{Smets_1993}
Smets Ph., \emph{Belief functions: The disjunctive rule of combination and the Generalized Bayesian Theorem}, Int. J. Approximate Reasoning, 9, pp. 1-35, 1993.
\bibitem{Zadeh_1978}
Zadeh L., \emph{Fuzzy sets as a basis for a theory of possibility}, Fuzzy sets and systems 1, pp. 3-28, 1978.
\end{thebibliography}
\end{document}